\documentclass[11pt]{amsart}

\newtheorem{Theorem}{Theorem}

\begin{document}

\title{Norm Bounds for Ehrhart Polynomial Roots}

\author{Benjamin Braun}
\address{Department of Mathematics, Washington University, St. Louis, MO}
\urladdr{math.wustl.edu/$\sim$bjbraun}
\email{bjbraun@math.wustl.edu}
\date{September 9, 2006.}

\begin{abstract}
M. Beck, J. De Loera, M. Develin, J. Pfeifle and R. Stanley found that the roots of the Ehrhart polynomial of a $d$-dimensional lattice polytope are bounded above in norm by $1+(d+1)!$. We provide an improved bound which is quadratic in $d$ and applies to a larger family of polynomials.
\end{abstract}
\maketitle

Let $P$ be a convex polytope in $R^n$ with vertices in $Z^n$ and affine span of dimension $d$; we will refer to such polytopes as \textit{lattice polytopes} and to elements of $Z^n$ as \textit{lattice points}.  A remarkable theorem due to E. Ehrhart, \cite{E}, is that the number of lattice points in the $t^{th}$ dilate of $P$, for non-negative integers $t$, is given by a polynomial in $t$ of degree $d$ called the \textit{Ehrhart polynomial} of $P$.  We denote this polynomial by $L_P(t)$, and let $\mathrm{Ehr}_P(x)=\sum_{t\geq 0}L_P(t)x^t$ denote its associated rational generating function.  For more information regarding Ehrhart theory, see \cite{BR}.

In \cite{BDDPS}, it was shown that for a lattice polytope $P$ of dimension $d$, the roots of $L_P(t)$ are bounded above in norm by $1+(d+1)!$.  However, the authors suggested that a bound that is polynomial in $d$ should exist and questioned whether this is a property of Ehrhart polynomials in particular or of a broader class of polynomials (see Remark 4.4 on page 26 of \cite{BDDPS}).  Our answer is the following:

\begin{Theorem} If $f$ is a non-zero polynomial of degree $d$ with real-valued, non-negative coefficients when expressed with respect to the polynomial basis $$B_d := \left\{ {t+d-j \choose d} : 0 \leq j \leq d \right\},$$  then all the roots of $f$ lie inside the disc with center $\frac{-1}{2}$ and radius $d(d-\frac{1}{2})$.
\end{Theorem}

The link between this situation and Ehrhart polynomials is that for a polynomial $f$ of degree $d$ over the complex numbers, there always exist complex values $h_j$ so that $$\frac{\sum_{j=0}^{d}h_jx^j}{(1-x)^{d+1}}=\sum_{t\geq 0}f(t)x^t.$$ \noindent As a result, $f$ can be expressed as  $$f(t)=\sum_{j=0}^{d}h_j{t+d-j \choose d}.$$  This is easily seen by expanding the rational function as a formal power series.  We then apply the following theorem, originally due to R. Stanley:

\begin{Theorem} (see \cite{S} and \cite{BR}.)  If $P$ is a $d$-dimensional lattice polytope with $$\mathrm{Ehr}_P(x)=\frac{\sum_{j=0}^{d}h_jx^j}{(1-x)^{d+1}},$$  then the $h_j$ are non-negative integers.
\end{Theorem}

\noindent Thus, our result applies to Ehrhart polynomials and more generally to Hilbert polynomials of certain Cohen-Macaulay modules (see \cite{BrH}, Corollary 4.1.10).  Theorem 1 is proved as follows.

\begin{proof} Let $d$ be a positive integer, let $D_d:=\{z:|z+\frac{1}{2}|\leq d(d-\frac{1}{2})\}$, and let $f$ be as given in the proposition.  It is enough to show that for any complex number $z$ not in $D_d$ there exists an open half-plane with zero on the boundary containing $B_d(z):=\{{z+d-j \choose d}:0\leq j\leq d\}$, since this implies that $f(z)$ is a non-trivial, non-negative linear combination of elements in a common open half-plane and is hence non-zero.

Each element of $B_d(z)$ is given by the product of $\frac{1}{d!}$ and $d$ consecutive members of $M:=\{(z+d), (z+d-1),\ldots,(z-d+2),(z-d+1)\}$.  The elements of $M$ are contained in a disk $D(z)$ of diameter $2d-1$ centered at $z+\frac{1}{2}$.  We claim that if $|z+\frac{1}{2}|>d(d-\frac{1}{2})$, which holds for $z\notin D_d$, then the angular width of $D(z)$ is less than $\frac{\pi}{d}$.  To see this, consider one of the lines through the origin tangent to $D(z)$.  The triangle formed by the origin, the point of tangency, and $z+\frac{1}{2}$ is a right triangle with hypotenuse of length $|z+\frac{1}{2}|$ and a side of length $d-\frac{1}{2}$ opposite the interior angle formed at the origin.  Hence, the interior angle at the origin is $\sin^{-1}\left(\frac{d-\frac{1}{2}}{|z+\frac{1}{2}|}\right)$, and thus the total angular width of $D(z)$ is $2\sin^{-1} \left(\frac{d-\frac{1}{2}}{|z+\frac{1}{2}|}\right)$.  Finally, we see that $$2\sin^{-1} \left(\frac{d-\frac{1}{2}}{|z+\frac{1}{2}|} \right) < 2\sin^{-1} \left(\frac{d-\frac{1}{2}}{d(d-\frac{1}{2})}\right) = 2\sin^{-1}\left(\frac{1}{d}\right) < \frac{\pi}{d}.$$

Therefore, the elements of $M$ all lie in a cone in the plane with apex the origin and angle width less than $\frac{\pi}{d}$.  Thus, the angular difference between $(z+d-j)\cdots(z-j+1)$ and $(z+d-j-1)\cdots(z-j)$ is less than $\frac{\pi}{d}$ for any $j$, $0\leq j < d$.  Hence, $B_d(z)$ lies in an open half-plane and our proof is complete.

\end{proof}

All the polynomials in $B_d$ have roots contained in $\{-d,-d+1,\ldots,d-1\}$.  For $1\leq j \leq d$, the number of polynomials in $B_d$ with $-j$ as a root is equal to the number with $-1+j$ as a root.  Thus, the location of the center of the disc in our proposition should not come as a surprise since the roots of the elements of $B_d$ are highly symmetric with respect to the point $\frac{-1}{2}$.  The line $x=\frac{-1}{2}$ also plays a prominent role for Ehrhart polynomials of cross-polytopes, as shown in \cite{BCKV} and \cite{R}.

It is interesting that our result only depends on $f$ having a ``nice'' representation with respect to $B_d$.  In our situation, the reason that $B_d$ is better than the standard monomial basis is that each of the polynomials in $B_d$ is of full degree $d$, and hence each such polynomial has $d$ roots.  In fact, by adapting our method one can obtain root bounds for any family of functions given by non-negative linear combinations of elements of a basis for degree $d$ polynomials that consists only of polynomials of degee $d$ having positive real leading coefficients and whose roots are known.

Thanks to John Shareshian for suggestions and advice, Matthias Beck and Sinai Robins for introducing me to Ehrhart theory, an anonymous referee for thoughtful comments, and Laura Braun for support and encouragement.


\begin{thebibliography}{10}
	\bibitem{BDDPS}{M. Beck, J. De Loera, M. Develin, J. Pfeifle, and R. Stanley, Coefficients and Roots of Ehrhart Polynomials, in Integer points in polyhedra --  geometry, number theory, algebra, optimization, volume 374 of Contemp. Math., pp. 15--36.  Amer. Math. Soc., Providence, RI, 2005.  arxiv:math.CO/0402148}
	\bibitem{BR}{M. Beck and S. Robins, Computing the Continuous Discretely, to be published by Springer books, draft available at math.sfsu.edu/beck/ccd.html.}
	\bibitem{BrH}{W. Brunz and J. Herzog, Cohen-Macaulay rings, Cambridge: Cambridge University Press, 1993.}
	\bibitem{BCKV}{D. Bump, K.-K Choi, P. Kurlberg, and J. Vaalar, A local Riemann hypothesis, I, Math. Z. 233 (2000), no. 1, 1--19.}
	\bibitem{E}{E. Ehrhart, Sur les poly\`edres rationnels homoth\'etiques \`a {$n$}\
               dimensions, C. R. Acad. Sci. Paris 254 (1962), 616--618.}
 	\bibitem{R}{F. Rodriguez-Villegas, On the zeros of certain polynomials, Proc. Amer. Math. Soc. 130 (2002), no. 8, 2251--2254.}
 	\bibitem{S}{R. Stanley, Decompositions of rational convex polytopes, Ann. Discrete Math. 6 (1980), 333--342.   Combinatorial mathematics, optimal designs and their applications (Proc. Sympos. Combin. Math. and Optimal Design, Colorado State Univ., Fort Collins, Colo., 1978).}
\end{thebibliography}
\end{document}